\documentclass[11pt,fullpage, doublespace]{amsart}
\usepackage{amsfonts, amssymb, amsmath}
\usepackage[all]{xy}
\usepackage{bm}
\usepackage{authblk}
\newtheorem{prop}{Proposition}[section]
\newtheorem{theorem}{Theorem}[section]

\newtheorem{defi}{Definition}[section]
\newtheorem{corr}{Corollary}[section]

\newcommand{\gm}{\gamma}
\newcommand{\dt}{\delta}
\newcommand{\ld}{\lambda}

\newcommand{\Z}{\mathbb{Z}}

\newcommand{\Q}{\mathbb{Q}}

\newcommand{\mapto}{\longrightarrow}

\title{On the Order of $a$ modulo $n$, on Average}
\author{Kim, Sungjin}
\affil[1]{Department of Mathematics \\ University of California, Los Angeles \\ Math Science Building 6617A \\ E-mail:  707107@gmail.com }

\makeatletter
\let\authors\AB@authors
\makeatother

\begin{document}
    \maketitle

    \begin{abstract}
    Let $a>1$ be an integer. Denote by $l_a(n)$  the multiplicative order of $a$ modulo integer $n\geq 1$.
    We  prove  that there is a positive constant $\dt$ such that if $x^{1-\dt} = o(y)$, then
    $$
    \frac1y\sum_{a<y}\frac1x\sum_{\substack{{a<n<x}\\{(a,n)=1}}}l_a(n) =   \frac x{\log x}\exp \left(B\frac{\log\log x}{\log\log\log x}(1+o(1))\right)$$
    where
    $$
    B=e^{-\gm}\prod_p \left(1-\frac 1{(p-1)^2(p+1)}\right).$$
    It was known for $y=x$ in ~\cite[Page 3]{KP} in which they refer to ~\cite{LS}.
        \end{abstract}

\section{Introduction}
    Let $a>1$ be an integer. If $n$ be coprime to $a$, we write  $d=l_a(n)$ if $d$ is the multiplicative order of $a$ modulo $n$. Then $d$ is the smallest positive integer in the congruence $a^d\equiv 1$ (mod $n$).

    The Carmichael's lambda function $\ld(n)$ is defined by the exponent of the group $(\Z/n\Z)^{*}$. It was known in ~\cite{EPS} that
    $$\frac1x \sum_{n<x} \lambda(n) = \frac x{\log x}\exp \left(B\frac{\log\log x}{\log\log\log x}(1+o(1))\right).$$
    Assuming GRH for Kummer extensions $\Q(\zeta_d, a^{1/d})$, P. Kurlberg and C. Pomerance ~\cite{KP} showed that
    $$
    \frac1x \sum_{n<x} l_a(n) = \frac x{\log x}\exp \left(B\frac{\log\log x}{\log\log\log x}(1+o(1))\right)$$
    with $B=e^{-\gm}\prod_p \left(1-\frac 1{(p-1)^2(p+1)}\right)$.
    The upper bound implicit is unconditional because $l_a(n)\leq \ld(n)$. An unconditional average result over all possible nonzero residue classes is obtained by F. Luca and I. Shparlinski ~\cite{LS}:
    $$
    \frac{1}{x} \sum_{n<x}\frac{1}{\phi(n)}\sum_{a<n} l_a(n) = \frac x{\log x}\exp \left(B\frac{\log\log x}{\log\log\log x}(1+o(1))\right).
    $$
    As pointed out in ~\cite{KP}, by partial summation, we have the following statistics on average order:
    $$\frac{1}{x^2}\sum_{a<x}\sum_{a<n<x}  l_a(n) =  \frac x{\log x}\exp \left(B\frac{\log\log x}{\log\log\log x}(1+o(1))\right).$$
    For fixed $a$, it seems that it is very difficult to remove GRH in P. Kurlberg and C. Pomerance's result with current knowledge. However, we expect that averaging over $a$ would give some information. So, we take average over $a<y$, but we do not want to have too large $y$ such as $y>x$. For all the average results in this paper, we assume that $y<x$, and try to obtain $y$ as small as possible.
    By applying a deep result on exponential sums by Bourgain ~\cite{B}, we prove the unconditional average result on a shorter interval.
    \begin{theorem}
    There is a positive constant $\dt$ such that, if $x^{1-\dt}=o(y)$, then
    $$
    \frac1y\sum_{a<y}\frac1x\sum_{\substack{{a<n<x}\\{(a,n)=1}}}l_a(n) =   \frac x{\log x}\exp \left(B\frac{\log\log x}{\log\log\log x}(1+o(1))\right)$$
    where
    $$
    B=e^{-\gm}\prod_p \left(1-\frac 1{(p-1)^2(p+1)}\right).$$
    \end{theorem}

   \section{Backgrounds}
    \subsection{Equidistribution}
    A sequence $\{a_n\}$ of real numbers are said to be equidistributed modulo $1$ if the following is satisfied:
    \begin{defi}
    Let $0\leq a<b\leq 1$. Suppose that
    $$
    \lim_{N\rightarrow\infty}\frac1N|\{n\leq N \ : \ a_n \in (a,b) \  \mathrm{mod}\  1\}| = b-a.$$
    Then we say that $\{a_n\}$ is equidistributed modulo $1$.
    \end{defi}
    A well-known criterion by Weyl ~\cite{W} is
    \begin{theorem}
    For any integer $k\neq 0$, suppose that
    $$
    \lim_{N\rightarrow\infty}\frac1N \sum_{n\leq N}e^{2\pi i ka_n } = 0.$$
    Then the sequence $\{a_n\}$ is equidistributed modulo $1$.
    \end{theorem}
    There was a series of efforts to obtain a quantitative form of the equidistribution theorem. Erd\H{o}s and Tur\'{a}n ~\cite{ET} succeeded in obtaining the following result:
    \begin{theorem}
    Let $\{a_n\}$ be a sequence of real numbers. Then for some positive constants $c_1$ and $c_2$,
    $$
    \sup_{0\leq a<b\leq 1} \left||\{ n\leq N :  a_n \in (a,b) \  \mathrm{mod} \ 1 \}| - (b-a)N\right| \leq c_1 \frac{N}{M+1} + c_2 \sum_{m=1}^M\frac1m \left|\sum_{n\leq N}e^{2\pi i m a_n}\right|.
    $$
    \end{theorem}
    H. Montgomery ~\cite{M} obtained $c_1=1$, $c_2=3$. C. Mauduit, J. Rivat, A. S\'{a}rk\H{o}zy ~\cite{MRS} obtained $c_1=c_2=1$. Thus, we have a quantitative upper bound of discrepancy when we have good upper bounds for exponential sums.
    \subsection{Exponential Sums in $\Z_n^{*}$}
     We define  arithmetic functions $a_n(d)$ and  $b_n(d)$ for $1\leq d|\ld(n)$ as follows:
    $$a_n(d)=|\{0<a<n: \ l_a(n)=d\}|,$$
    $$b_n(d)= |\{0<a<n: \ a^d\equiv 1 \ (\mathrm{mod} \ n)\}|.$$
    Then $$a_n(d)=\sum_{d'|d} \mu\left(\frac d{d'}\right) b_n(d').$$
    We give some algebraic remarks about the function $b_n(d)$. First, we see that
    $$H_{n,d}:=\{0<a<n: \ a^d\equiv 1 \ (\mathrm{mod} \ n)\}$$ forms a subgroup of $\Z_n^{*}$ of order $b_n(d)$. The following proposition is from elementary group theory:
    \begin{prop}
    Let $H_{n,d}$ and $b_n(d)$ be defined as above. For any $k|n$, denote by $\pi_k$ the reduction modulo $n/k$. Then we have
    $$
    \pi_k: H_{n,d} \mapto H_{n/k, d}
    $$
    where $\pi_k$ is a group homomorphism with kernel
    $$
    K = \{ 0<a<n: \ a^d\equiv 1 (n), \ a\equiv 1 (n/k) \}.
    $$
    By the First Isomorphism Theorem, we have
    $$
    |K| = \frac{b_n(d)}{|\pi_k(H_{n,d})|}\leq k.$$
    \end{prop}
    Note that the map $\pi_k$ restricted to $H_{n,d}$ is not always surjective. To see this, let $p>2$ prime number, and $a=p+1$, $d=p$, $n/k=p^2$, $n=p^3$. Then
    $$
    a^p \equiv p^2 + 1 \ (\mathrm{ mod }\ p^3).$$
    Thus,
    $$
    a^p \equiv 1 \ (\mathrm{ mod }\ p^2).$$
    But for any $a'\equiv  a \ (\mathrm{mod} \ p^2)$, so that $a'=p^2j + p+1$ for some integer $j$, we have
    $$
    (a')^p \equiv (p+1)^p  \ (\mathrm{ mod }\ p^3) \equiv p^2 + 1 \ (\mathrm{mod} \ p^3).$$
    From this, we see that the element $a=p+1\in H_{n/k}$ is not a preimage of $\pi_k$.
The proof of $|K|\leq k$ is clear by $a\equiv 1 (n/k)$.

    J. Bourgain ~\cite{B} proved a nontrivial exponential sum result when a subgroup $H$ of $\Z_n^{*}$ has order greater than $n^{\epsilon}$ for  $\epsilon>0$.
    \begin{theorem}
    Let $n\geq 1$. For any $\epsilon>0$, there exist  a constant $\dt=\dt(\epsilon)>0$ such that for any subgroup $H$ of $\Z_n^{*}$ with $|H|>n^{\epsilon}$,
    $$
    \max_{(m,n)=1}\left|\sum_{a\in H}e^{2\pi i m\frac {a}n}\right|<n^{-\dt}|H|.$$
    \end{theorem}
    \begin{corr}
    Let $\epsilon>0$ be arbitrary, and let $y\geq 1$. Assume that $d|\ld(n)$ and $b_n(d) > n^{\epsilon}$. Then there exists $\dt=\dt(\epsilon)>0$ such that
    $$
    \sum_{a<y, \ a^d\equiv 1 (n)} 1 = \frac yn b_n(d) + O(b_n(d)n^{-\dt}).
    $$
    \end{corr}
    If $d|\ld(n)$, the congruence $a^d\equiv 1$ yields $b_n(d)$ roots in $\Z_n$. Thus, we need to count $a<y$ satisfying those $b_n(d)$ congruences modulo $n$. Considering $\frac yn = \lfloor \frac yn \rfloor + \frac yn - \lfloor \frac yn \rfloor$, it is enough to prove the result for $y<n$. We apply the Erd\H{o}s-Tur\'{a}n inequality to the set $\{\frac an : 0<a<n, \  a^d \equiv 1 (n)\}$. Then
    \begin{align*}
    \left||\{ 0<a<n :  a^d\equiv 1 (n), \frac an \in (0,\frac yn) \textrm{ mod $1$} \}| - \frac yn b_n(d)\right| &\leq   \frac{b_n(d)}{n} +   \sum_{m=1}^{n-1}\frac{1}m \left|\sum_{a \in \Z_n , \ a^d\equiv 1 (n)}e^{2\pi i m \frac an}\right|.
    \end{align*}
    Unlike the prime modulus case, we immediately encounter a problem. The exponential sum result (Theorem 2.3) is only for $(m,n)=1$, but the sum takes all $1\leq m<n$. Then we have too many terms with $(m,n)\neq 1$. Therefore, we need some modification in applying the Erd\H{o}s-Tur\'{a}n inequality. A starting point is to observe that we can take $M$ arbitrary in the Erd\H{o}s-Tur\'{a}n inequality.

\it Proof of Corollary 2.1) \rm

    Assuming that $k|n$ and $b_n(d)>n^{\epsilon}$, we have
    \begin{align*}
    n^{\epsilon} < b_n(d) \leq k |\pi_k(H_{n,d})|.
    \end{align*}
    Then
    \begin{align*}
    \frac{n^{\epsilon}}k<|\pi_k(H_{n,d})|.
    \end{align*}
    If we can assume that
    $$
    \left(\frac nk\right)^{\epsilon''} < \frac{n^{\epsilon}}k$$
    for some positive $\epsilon'' < \epsilon$, then we can use Theorem 2.3 with $\epsilon''$ and $\dt''=\dt(\epsilon'')$. This is achieved by
    $$
    k<n^{\frac{\epsilon -\epsilon''}{1-\epsilon''}}.$$
    Let $\epsilon'=\frac{\epsilon -\epsilon''}{1-\epsilon''}$ and we take  $M+1=\lfloor n^{\epsilon'}\rfloor$ in the Erd\H{o}s-Tur\'{a}n inequality. Then we have reduced the number of terms appearing in the sum on the right side. We rewrite the sum by substituting $(m,n)=k$, $\frac mk = j$ and apply Theorem 2.3 to the exponential sums inside. This is possible due to
    $$
    \left(\frac nk\right)^{\epsilon''} < |\pi_k(H_{n,d})|
    $$
    and $\pi_k(H_{n,d})$ being a subgroup of $\Z_{n/k}^{*}$.
    The sum on the right becomes
    \begin{align*}
    \sum_{m<n^{\epsilon'}}\frac1m \left|\sum_{a\in H_{n,d}}e^{2\pi i m \frac an}\right|&\leq\sum_{\substack{{k|n}\\{k<n^{\epsilon'}}}}\frac1{k }\sum_{\left(j,\frac nk\right)=1}\frac1j\left|\sum_{a\in H_{n,d}}e^{2\pi i j \frac a{n/k}} \right|\\
    &=\sum_{\substack{{k|n}\\{k<n^{\epsilon'}}}}\frac1{k }\sum_{\left(j,\frac nk\right)=1}\frac1j\frac{b_n(d)}{|\pi_k(H_{n,d})|}\left|\sum_{a\in\pi_k(H_{n,d})}e^{2\pi i j \frac a{n/k}} \right|\\
    &\leq \sum_{\substack{{k|n}\\{k<n^{\epsilon'}}}}\frac1{k }\sum_{\left(j,\frac nk\right)=1}\frac1j\frac{b_n(d)}{|\pi_k(H_{n,d})|}|\pi_k(H_{n,d})| \left(\frac nk\right)^{-\dt''}\\
    &\leq n^{-\dt''(1-\epsilon')}b_n(d)(1+\log n)^2.
    \end{align*}
    Thus, the Erd\H{o}s-Tur\'{a}n inequality gives
    $$
    \left||\{ 0<a<n :  a^d\equiv 1 (n), \frac an \in (0,\frac yn) \textrm{ mod $1$} \}| - \frac yn b_n(d)\right|\leq \frac{b_n(d)}{n^{\epsilon'}}+ b_n(d)n^{-\dt''(1-\epsilon')}(1+\log n)^2.
    $$
    Therefore we can take $0< \dt < \textrm{min} (\epsilon', \dt''(1-\epsilon'))$. This completes the proof of Corollary 2.1.

    Corollary 2.1  plays  a key role  in proving Theorem 1.1. Note that the upper bound provided in Corollary 2.1  is significantly better than the trivial bound which is:
    $$
    \sum_{a<y, \ a^d\equiv 1 (n)} 1 = \frac yn b_n(d) + O(b_n(d)).
    $$
    \section{Proof of Theorems}

    \subsection{Proof of Theorem 1.1}
    We start with the change of order in summation:
    \begin{align*}
    \sum_{a<y}\sum_{n<x}l_a(n)   &=\sum_{d<x}d \sum_{\substack{{n<x}\\{d|\ld(n)}}}\sum_{\substack{{a<y}\\{l_a(n)=d}}} 1\\
    &=\sum_{d<x}  d \sum_{\substack{{n<x}\\{d|\ld(n)}}}\sum_{\substack{{d'|d}\\{b_n(d')<n^{\epsilon}}}}\mu\left(\frac{d}{d'}\right)\sum_{\substack{{a<y}\\{a^{d'}\equiv 1 (n)}}}1 + \sum_{ d<x}d \sum_{\substack{{n<x}\\{d|\ld(n)}}}\sum_{\substack{{d'|d}\\{b_n(d')\geq n^{\epsilon}}}}\mu\left(\frac{d}{d'}\right)\sum_{\substack{{a<y}\\{a^{d'}\equiv 1 (n)}}}1 \\
    &=\sum_{ d<x}d \sum_{\substack{{n<x}\\{d|\ld(n)}}} \sum_{\substack{{d'|d}\\{b_n(d')< n^{\epsilon}}}}\mu\left(\frac{d}{d'}\right)\left(\frac yn b_n(d') + O(b_n(d'))\right)\\
    & \ \ + \sum_{ d<x}d \sum_{\substack{{n<x}\\{d|\ld(n)}}} \sum_{\substack{{d'|d}\\{b_n(d')\geq n^{\epsilon}}}}\mu\left(\frac{d}{d'}\right)\left(\frac yn b_n(d') + O(b_n(d')n^{-\dt})\right)\\
    &=\sum_{ d<x}d \sum_{\substack{{n<x}\\{d|\ld(n)}}} \frac yn a_n(d) + O(E_1)+O(E_2),
    \end{align*}
    where
    \begin{align*}
    E_1 &= \sum_{d<x}d\sum_{\substack{{n<x}\\{d|\ld(n)}}}\sum_{\substack{{d'|d}\\{b_n(d')<n^{\epsilon}}}}\left|\mu\left(\frac{d}{d'}\right)\right|b_n(d')\\
    &\ll \sum_{d<x}d\sum_{\substack{{n<x}\\{d|\ld(n)}}} \sum_{d'|d} n^{\epsilon}\\
    &\ll \sum_{d<x}d\tau(d) \sum_{\substack{{n<x}\\{d|\ld(n)}}}n^{\epsilon}\\
    &= \sum_{n<x}n^{\epsilon}\sum_{d|\ld(n)}d\tau(d)\\
    &\ll x^{2+\epsilon+o(1)}  ,
    \end{align*}
    and
    \begin{align*}
    E_2 &= \sum_{d<x}d\sum_{\substack{{n<x}\\{d|\ld(n)}}}\sum_{\substack{{d'|d}\\{b_n(d')\geq n^{\epsilon}}}}\left|\mu\left(\frac{d}{d'}\right)\right|b_n(d')n^{-\dt}\\
    &\ll \sum_{d<x}d\sum_{\substack{{n<x}\\{d|\ld(n)}}} b_n(d) n^{-\dt} \sum_{d'|d} 1\\
    &= \sum_{n<x} \sum_{d|\ld(n)}d b_n(d) \tau(d)n^{-\dt}\\
    &\le \sum_{n<x} n^{1-\dt} \sum_{d|\ld(n)}d\tau(d)\\
    &\ll x^{3-\dt+o(1)}.
    \end{align*}
    Now we treat the main term:
    \begin{align*}
    \sum_{  d<x}d \sum_{\substack{{n<x}\\{d|\ld(n)}}} \frac 1n a_n(d)
    &=\sum_{n<x} \frac{1}n\sum_{d|\ld(n)} da_n(d).
    \end{align*}
    Taking $\dt$ to satisfy $2+\epsilon\le 3-\dt$, we have
    $$
    \sum_{a<y}\sum_{n<x} l_a(n)  = y \sum_{n<x} \frac{1}n\sum_{d|\ld(n)} da_n(d) + O(x^{3-\dt+o(1)} ).
    $$
    Let $u(n)=\frac1{\phi(n)} \sum_{d|\ld(n)} da_n(d)$ be the average multiplicative order of the elements of $(\Z/n\Z)^{*}$. The following is proven in ~\cite[Theorem 6]{LS}:
    \begin{theorem}
    $$\frac1x \sum_{n<x} u(n)=\frac x{\log x}\exp \left(B\frac{\log\log x}{\log\log\log x}(1+o(1))\right).$$
    \end{theorem}
    What we have for the main term is the middle term in the following inequalities:
    $$
    \frac1{\log\log x}\sum_{n<x} u(n) \ll \sum_{n<x} \frac{\phi(n)}n u(n) \leq \sum_{n<x} u(n).$$
    Since $\log\log\log x = o\left( \frac{\log\log x}{\log\log\log x}\right)$, it follows that
    $$\sum_{n<x} \frac{\phi(n)}n u(n)=\frac {x^2}{\log x}\exp \left(B\frac{\log\log x}{\log\log\log x}(1+o(1))\right).$$
    Hence, we have
    $$
    \sum_{a<y}\sum_{n<x} l_a(n) = \frac {yx^2}{\log x}\exp \left(B\frac{\log\log x}{\log\log\log x}(1+o(1))\right) + O(x^{3-\dt+o(1)} ).
    $$
    Moreover, if for some $0<\dt'<\dt$, and $x^{1-\dt'} = o(y)$, then the error term can be included in the term with $o(1)$.  The terms that appear when $n\leq a$, are also included in the term with $o(1)$. This completes the proof of Theorem 1.1.
        \flushleft

\end{document}